\input amstex
\input xy
\xyoption{all}
\documentstyle{amsppt}
\document
\magnification=1200
\NoBlackBoxes
\nologo
\hoffset1.5cm
\voffset2cm
\vsize16cm

\bigskip

\centerline{\bf  HIGHER STRUCTURES,}

\medskip

\centerline{\bf  QUANTUM GROUPS,}

\medskip

\centerline{\bf AND GENUS ZERO MODULAR OPERAD}

\bigskip

\centerline{\bf Yuri I.~Manin}

\medskip

\centerline{Max--Planck--Institut f\"ur Mathematik, Bonn, Germany}

\bigskip

\hfill{\it To Fedor}
\medskip
{\bf Abstract.} In my Montreal lecture notes of 1988, it was suggested that the theory of linear quantum groups
can be  presented in the framework of the category of {\it quadratic algebras} (imagined as
algebras of functions on ``quantum linear spaces"), and  quadratic algebras
of their inner (co)homomorphisms.

Soon it was understood (E.~Getzler and J.~Jones, V.~ Ginzburg, M.~Kapranov, M.~Kontsevich,
M.~Markl, B.~Vallette et al.)
that the class of {\it quadratic operads} can be introduced
and the main theorems about quadratic algebras can be generalised to the
level of such operads, if their components  are {\it linear spaces} (or objects of more general monoidal categories.)

When quantum cohomology entered the scene,
it turned out that  the basic tree level (genus zero)
(co)operad of quantum cohomology not only is {\it quadratic} one, but its
{\it components are themselves quadratic algebras.}

In this short note, I am studying the interaction of quadratic algebras structure with operadic
structure
in the context  of enriched 
category formalism due to G.~M.~Kelly et al.

\bigskip


\medskip{


\bigskip

\centerline{\bf  1. Motivation and brief summary}

\medskip

{\bf 1.1. Motivation.} In the mathematical language of quantum physics, an important role is played by a list of 
different quantisation formalisms prescribing how to pass from a classical
description of an ``isolated system'' to its quantum description.

\smallskip
In all items of this list, from Schr\"odinger's commutation relations
to Feynman's path integrals, the Planck constant (actually, not a number but a natural 
unit of measurement for action) appears as a ``small deformation
parameter'' of the relevant classical structure.
\smallskip

When it was recognised that {\it symmetry groups} of classical laws
can and must be quantised as well (see [Ji85], [Dr86], [FadReTa87]),
the modern theory of {\it quantum groups} emerged. 
In the Montreal lecture notes [Ma88], I suggested that the theory of linear quantum groups
can be  presented in the framework of the category of quadratic algebras imagined as
algebras of functions on "quantum linear spaces", and their internal
{\it hom--} and {\it cohom}--objects. See also [Ma87] and [Ma91].

\smallskip

Similarly, from the great insights of quantum filed theory regarding
string models, there emerged various chapters of mathematical theory of {\it quantum
(co)homology} where a symplectic or algebraic manifold $M$ is studied
by considering maps of Riemannian surfaces/algebraic curves to $M$
and deforming 
topological/motivic invariants of the manifold using all such maps.

\smallskip

Mathematical formalism of quantum cohomology involves 
the construction of an operad whose components are
(co)homology groups of moduli spaces $\overline{M}_{g,n}$ of stable pointed curves.

\smallskip

This stimulates a search of  a  ``noncommutative", or ``quantum'' version of
the operad of cohomology/motives $\{ h(\overline{M}_{g,n})\}$ of moduli spaces: 
cf. an attempt in [DoShVa15] (for $g=0.$)

\smallskip

Actually, a considerable part of this
generic formalism is already implicit in the vast framework of ``generalised operads''
developed recently: see [KaWa17] and an earlier project [BorMa07].
For other aspects of categorical quantisation, cf. [DaSt04].
\smallskip

In this short note, I focus on the genus zero, or {\it tree level} part
of the big modular operad. The point is that 
 all its components are not just linear spaces but quadratic rings,
and quadraticity furnishes many additional beautiful structures that
are absent in more general contexts.

Moreover, this operad itself is a {\it quadratic} one, and
properties of such operads   generalise properties
of quadratic algebras such as existence of Koszul duality constructions
(cf.  [GiKa94], [Va08], and the monograph [LoVa12]).

\medskip

{\bf 1.2. Summary.} Section 2 is a brief survey of main definitions and constructions
of the theory of quadratic algebras and their category $\bold{QA}.$

\smallskip

Section 3 introduces the first higher structure related to quadraticity: 
$\bold{QA}$--{\it enrichment of} $\bold{QA}$ in the style of Kelly [Ke82].

\smallskip

The main statement of section 4 is a description of of  quantised version
of the tree level modular (co)operad, where quantisation here is understood
as lifting the  (co)operadic (co)multiplications to the level of Kelly's 
enrichment of $\bold{QA}$.

\smallskip

For higher structures of different types , see [DoShVa12].

\smallskip

This note was motivated  by reflections on the very badly understood
``self--referentiality'' of the genus zero modular operad: cf. [MaSm13,14].
At the last subsection of sec. 4, I briefly comment on this problem.

\medskip

{\it Acknowledgements.} I am very grateful to Bruno Vallette, Ralph Kaufmann, 
Vladimir Dotsenko, and Mikhail Kapranov who read
preliminary versions of this note and greatly helped me  clarify and/or correct
some misguided statements.

\bigskip

\centerline{\bf 2. Category of quadratic algebras}

\medskip

{\bf 2.1 General notations.} Below $K$ denotes a fixed (commutative) field of characteristic zero. We start with
a category of vector spaces over $K$. 
\smallskip

{\it A quadratic algebra}  is a graded $K$--algebra $A=\oplus_{i=0}^{\infty} A_i$,
where $A_0=K$, $A_1$ is a finite dimensional subspace generating $A$,
and such that an appropriate subspace $R(A)\subset A_1^{\otimes 2}$
generates the ideal of all relations between elements of $A_1$.
In other words, $A$ is given together with the surjective morphism of the tensor algebra
of $A_1$ to $A$,
  whose kernel in the component of degree $d\ge 2$ equals 
  $$
  \sum_{i+k=d-2} A_1^{\otimes i}\otimes_K R(A) \otimes_K   A_1^{\otimes k}.
  $$
\smallskip
It is often convenient to write, as in [Ma88] $A\leftrightarrow (A_1, R(A))$.
\smallskip
Quadratic algebras are objects of the category $\bold{QA}$,
in which morphisms $A\to B$ can be described as linear maps 
$f:\,A_1 \to B_1$ such that $(f\otimes f) (R(A)) \subset R(B).$

\smallskip

There is also the natural functor  $\bold{QA} \to \bold{Lin}_K$ (where $\bold{Lin}_K$
is the category of finite dimensional
linear spaces over $K$). It is given by $A\mapsto A_1$.
\smallskip

Whenever we are dealing only with $\bold{QA}$ as a category, we may simply
denote its objects $(A_1, R(A))$.
\smallskip

The key difference between  $\bold{Lin}_K$ and $\bold{QA}$ can be briefly summarised
as follows: the former category has {\it two} natural symmetric monoidal
structures $\otimes_K$ and $\oplus$, whereas $\bold{QA}$ has at least {\it four}
such structures $\bullet ,\circ , \otimes , \underline{\otimes}$. These structures
are obtained from $\bold{Lin}_K$  by liftings via the contravariant {\it duality functor} $!$\,
which is defined on objects by
$$
(A_1, R(A))^! := (A_1^*, R(A)^{\bot} ).
$$
Here ${}^{*}$ means the linear duality functor in $\bold{Lin}_K$, and $\bot$ denotes
the appropriate orthogonal complement functor.

\smallskip

For precise definitions and discussion of their basic properties needed below, we
will mostly refer to [Ma88], especially sec. 3 and 4, pp.~19--28.
See also [Ma87], [Ma17] and [BoyDr13].

\smallskip

Whenever convenient, we will work with strict monoidal categories
and write $=$ for canonical identifications. 
\bigskip

\centerline{\bf 3. $\bold{QA}$--enrichment of $\bold{QA}$}

\medskip

{\bf 3.1. Monoidal categories.} We will use the definition of (symmetric) {\it  monoidal category}
as it was stated in [Ke82] (online version, sec.~1.1).
Similar, or equivalent, notions were also considered in
[DeMi82] under the name of {\it tensor category},
and in many other papers. 
\smallskip

Briefly, a monoidal category is a category $\Cal{V}_0$, endowed with 
a bifunctor ``tensor product'' $\otimes:\,  \Cal{V}_0 \times \Cal{V}_0 \to\Cal{V}_0$,
satisfying compatible associativity and commutativity constraints
(in [Ke82], they are coherence axioms (1.1), (1.14), (1.15)). 
Moreover, the structure data of a monoidal category
must include a {\it unit object} $I$ together with diagrammatic
expressions of the fact, that it is the left and the right identity
for the tensor product: see [Ke82], (1.16).

\smallskip

An additional condition in the treatment
of enrichment by monoidal categories in [Ke82]  is the idea of its {\it closedness}.
A monoidal category is called {\it closed} if each functor of right tensor multiplication by a fixed object
$* \mapsto *\otimes Y$ has a right adjoint 
$* \mapsto [Y,*]$, that is ([Ke82], (1.23)):
$$
Hom_{\Cal{V}_0}(X\otimes Y, Z) = Hom_{\Cal{V}_0} (X, [Y,Z]).
$$

Kelly also  introduces  unit and counit functors
$$
d:\, X\mapsto [Y, X\otimes Y],\quad\quad  e: [Y,Z]\otimes Y \mapsto [Y,Z]\otimes Z :
$$
see [Ke82], (1.24).

\smallskip

\smallskip

Below, we will consider the symmetric monoidal category 
$(\bold{QA}, \bullet )$ with unit $K[t]/(t^2)$, where the {\it black product} $\bullet$
is defined on objects by
$$
(A_1, R(A))\bullet (B_1, R(B)):= (A_1\otimes_K B_1, S_{(23)}(R(A)\otimes_K R(B) ))
$$

We will need also the $!$--dual {\it white product} $\circ$ defined on objects by
$$
(A_1, R(A))\circ (B_1, R(B)):= (A_1\otimes_K B_1, S_{(23)}(R(A) \otimes_K B_1^{\otimes 2} +
A_1^{\otimes 2} \otimes_K R(B) )).
$$

\medskip

{\bf 3.2. Proposition.} {\it a) In the monoidal category $(\bold{QA}, \bullet )$
we have  functorial identifications
$$
Hom_{\bold{QA}}(A \bullet B, C)   = Hom_{\bold{QA}}(A, B^! \circ C)
$$
 so that $[B,C]=   B^! \circ C$.
 
 \medskip
 
 b) This adjoint functor is endowed with  counit  $e$.
}

\medskip

{\bf Proof.} a) The adjointness property is proved in [Ma88] (sec.~4, Theorem 2, pp.~25--26).
More precisely, the morphism in $\bold{QA}$ induced by a  linear map $f:\,A_1\otimes B_1\to C_1$
can be identified with morphism in $\bold{QA}$ induced by the  linear map
 $g:\,A_1\to  B_1^{*}\otimes C_1$ as is standard in the category of vector spaces,
 and the compatibility with quadratic relations is checked in [Ma88] directly.
\smallskip

b) The quadratic algebra $[B,C]=   B^! \circ C$ is denoted in [Ma88], p.~26, by
$\underline{Hom} (B, C)$. The morphism $\beta$ on p.~27, and its properties
discussed there show that it has the respective counit.

\medskip

{\bf 3.3. The enrichment.} Here I will describe the $(\bold{QA},\bullet )$--category in the sense of [Ke82], sec.~1.2,
whose set of objects is the same as in $\bold{QA}$.

\smallskip

For any two quadratic algebras $A,B$, the respective {\it hom--object} of [Ke82]
will be $A^!\circ B$ denoted  $\underline{Hom} (A,B)$ in [Ma88], sec.~4.

\smallskip

The composition law (Kelly's $M_{ABC}$) is our morphism $\mu = \mu_{ABC}$
in [Ma88], (7), p.~26.
\smallskip
Kelly's  identity elements $j_A:\, K[t]/(t^2) \to A^!\circ A$ and their generalisations are
introduced and discussed in sec.~9 of [Ma88].

\smallskip

Finally, we must check the associativity and unit axioms for this enrichment 
expressed by Kelly's commutative diagrams (1.3) and (1.4).
They are all checked in [Ma88].

\bigskip

\def\P{\Cal{P}}
\def\QA{\bold{QA}}

\centerline{ \bf 4. Genus zero modular operad and its enrichment}

\medskip

{\bf 4.1. Operads and cooperads  in the category of quadratic algebras.} 
The notion of operad that we have in mind in this introductory
subsection is a special case of  constructions described in [BorMa07].
Briefly, an operad $\P$ is a tensor functor between symmetric monoidal
categories $(\Gamma ,\coprod )\to (\bold{QA}, \otimes )$
where $\Gamma$ is a category of labelled (finite) graphs with disjoint union $\coprod$,
and the tensor product in  $\QA$ is defined on objects  by
$$
(A_1, R(A)) \otimes (B_1, R(B)) := (A_1\oplus B_1, R(A)\oplus [A_1,B_1]\oplus R(B)).
$$
%Narrowing our definition even more, we will be considering {\it shuffle operads}
%in the sense of [LoVa12], Sec 8.2. 
\smallskip
The data completely determining such an operad is  the set of morphisms in the target category  
$(\bold{QA}, \otimes )$
$$
 \P(k)\otimes \P(m_1)\otimes \P(m_2)\otimes\cdots \otimes \P(m_k) \to
 \P(n), \ n= m_1+m_2+\dots +m_k
 \eqno(4.1)  
$$
indexed by unshuffles of $\{1,2,\dots n\}$. They are called
{\it operadic multiplications.}
For a list of axioms for them, cf.~ [LoVa12],
defining {\it classical operads},
see [LoVa12], Prop~5.3.1.
\smallskip

The relevant notion of cooperad is obtained by inversion of arrows in (4.1).
\medskip

{\bf 4.1.1. Operadic data in the $\QA$--enrichment of $\QA$.} According to the Proposition 3.2
above, for any three quadratic algebras $A,B,C$ we have a canonical identification
$$
Hom_{\QA}(A\bullet B, C)  =  Hom_{\QA}(A, B^! \circ C).
\eqno(4.2)
$$
Putting here $A=K[t]/(t^2)$ which is the unit object in $(\QA ,\bullet )$, we get
$$
Hom_{\QA}(B, C)  =  Hom_{\QA}(K[t]/(t^2), B^! \circ C)
$$
$$
= \{ d\in B_1^*\otimes C_1\,|\, S_{(23)}(d^{\otimes 2}) \in R(B)^{\bot}\otimes C_1^{\otimes 2} +
(B_1^*)^{\otimes 2} \otimes R(C)\} .
\eqno(4.4)
$$
This finally means that the Kelly enrichment of such an operad $\P$ is given by the family of quadratic algebras
$$
 (\P(k)\otimes \P(m_1)\otimes \P(m_2)\otimes\cdots \otimes \P(m_k))^!\circ
 \P(m_1+m_2+\dots +m_k) 
\eqno(4.5)  
$$
endowed with a family of elements in the linear spaces
$$
 (\P(k)\otimes \P(m_1)\otimes \P(m_2)\otimes\cdots \otimes \P(m_k))_1^*\otimes
 \P(m_1+m_2+\dots +m_k)_1 
 \eqno(4.6) 
$$
indexed by unshuffles.

\smallskip

Moreover, composition and associavity axioms for these elements can be also lifted
to the Kelly enrichment.

\medskip

{\bf 4.1.2. Cooperads and coproducts.} As Bruno Vallette pointed out to me, 
in various definitions of operads/cooperads in monoidal categories
(cf.~ [MarShSt02], Part II, Ch.~1)  it is often useful to require
the target  category to be endowed with
a coproduct. Moreover, monoidal product must be distributive 
over coproduct: see [LoVa12], sec. 5.3.5.

\smallskip

{\bf 4.1.3. Lemma.}   {\it In the monoidal category $(\bold{QA}, \bullet )$, the initial object is
$K$, and the coproduct $\times$ is given on objects by
 $$
 (A_1,R(A)) \times  (B_1,R(B)) := (A_1\oplus B_1, R(A)\oplus R(B)).
$$
}
The distributivity of $\bullet$ can be checked then by a short direct computation.

\smallskip

In the category of graded associative algebras, this coproduct is known as 
free product.

\medskip

{\bf  4.2. Genus zero modular operad.} In this subsection, I will describe the main 
motivating example  of the shuffle operad in the category $\QA$:  the genus zero modular (co)operad
(also called tree--level cyclic CohFT (co)operad) $P$.

\smallskip

 The component of arity $n$ for $n\ge 2$ of  $P$
is the cohomology ring $P(n):=H^*(\overline{M}_{0, n+1}, \bold{Q})$ where $\overline{M}_{0, n+1}$
is the moduli space (projective manifold) parametrising stable curves of genus zero
with $n+1$ labelled points. Component of arity 1 is $\bold{Q}$.

\smallskip
 Structure morphisms (cooperadic comultiplications)
$$
P(m_1+m_2+\dots +m_k)\to P(k)\otimes P(m_1)\otimes P(m_2)\otimes\cdots \otimes P(m_k)  
\eqno(4.7)  
$$
are maps induced by the maps of moduli spaces defined point--wise  by a glueing
of the respective stable curves:
$$
\overline{M}_{0, k+1}\times      \overline{M}_{0, m_1+1}\times \cdots
\times   \overline{M}_{0, m_k+1} \to   \overline{M}_{0, m_1+ \cdots +m_k+1}.
\eqno(4.8)
$$

\medskip

{\bf 4.2.1. Proposition.}  {\it a) For every $n\ge 3$,  $P(n)$ is a quadratic algebra
with linear space of generators $P(n)_1 = H^2(\overline{M}_{0,n+1})$ of
dimension 
$$
2^n - \frac{n(n+1)}{2}-1.
$$

b) Comultiplications (4.7) are morphisms of quadratic algebras.
}
\smallskip

{\bf Proof.} For a) and further details, see [Ma99], Ch.~III, sec.~3, and references therein.

\smallskip

Part b) follows from the Proposition 8.3 (a) of [Fu84]:
any morphism of smooth projective
manifolds $X\to Y$ induces a functorial  homomorphism of Chow rings
$f^*: A^*(Y)\to A^*(X)$. Indeed, $P(n):=H^*(\overline{M}_{0, n+1}, \bold{Q})$
are just Chow rings graded by algebraic codimension of respective cycles.

\medskip

{\bf 4.2.2. Remark.} Algebras classified/encoded by $P$, are directly described below,
in the Def.~4.4.1.

\smallskip

There is another interesting operad $G$ whose components of every arity
are quadratic algebras as well. It encodes {\it Gerstenhaber algebras:}
cf.~[LoVa12], pp.~506 and 536. Each $G(n)$ can be represented as 
the homology ring of the Fulton--MacPherson
compactification of the space of configurations of $n$ points in $\bold{R}^2.$

\medskip

{\bf 4.3. Additional information about $P$ and $P$--algebras.} I will briefly repeat here
a description of $P$  as it was given
 in [KoMa94] and later considerably  generalised in [BarMa07] and [KaWa17].
 
 \smallskip
 
 I start with the combinatorial definition of relevant graphs.

\smallskip
(i) A stable tree $\tau$ is a diagram of pairwise disjoint finite sets $(V_{\tau}, E_{\tau}, T_{\tau})$
and boundary maps 
$$
b_T:\, T_{\tau} \to V_{\tau}, \quad b_E:\, E_{\tau} \to \{unordered\ pairs\
of\ distinct\ vertices\}.
$$ 

\smallskip

A geometric realization of $\tau$ is a CW--complex whose 1--simplexes are (bijective to)
 $E_{\tau}\cup T_{\tau}$ ({\it edges and tails}) and $0$--simplexes are (bijective to) $V_{\tau}$
( {\it vertices.})
The geometric realisation of $\tau$ must be connected and simply--connected, i.~e.
to be a tree. Each vertex  must belong to the boundary of either one tail, or one tail and
$\ge 2$ edges, or else
or $\ge 3$ edges (stability condition).

\smallskip

(ii) Stable trees are objects of a category, in which every morphism $f:\,\tau \to \sigma$
consists of three maps
$$
f_v:\, V_{\tau} \to V_{\sigma},\quad 
f^t:\, T_{\sigma} \to T_{\tau},\quad
f^e:\, E_{\sigma} \to  E_{\tau}.
$$
satisfying conditions spelled out in [KoMa94], Definition 6.6.2.

\smallskip

(iii) Let now $F$ be a finite set of cardinality $\ge 3.$ Below we will denote by
$\overline{M}_{0,F}$ the moduli space of stable curves of arithmetic genus zero
endowed with a collection of pairwise different smooth points labelled by $F$.

\smallskip

One can define a functor $\Cal{M}$ from the category of stable trees above to the category
of projective algebraic manifolds. On objects, it is defined by
$$
\Cal{M}:\, \tau \mapsto \prod_{v\in V_{\tau}} \overline{M}_{0,F_{\tau}(v)}.
$$
 Here  $F_{\tau}$ denotes the set of flags of $\tau$ that is,  (pairs $\{ edge,\ one\ vertex\ of\ it\}$), and
 $F_{\tau}(v)$ denotes the set  of all flags, containing the vertex $v$.
 
 \smallskip
 
 For the definition of $\Cal{M}$ on morphisms, see [KoMa94], p.~555.
  
\medskip

{\bf 4.3.1. Operadic endomorphisms and $\Cal{M}$--algebras.} Let $L$ be an object 
of the category $ \bold{Lin}^s_K$ of
finite--dimensional $K$--linear superspaces with a non--degenerate even
scalar product. One can define the operad $\bold{OpEnd}\,L$ as the functor on stable trees
defined on objects by
$$
\bold{OpEnd}\,L\,(\tau ):= L^{\otimes F_{\tau}}.
$$
For the definition on morphisms, see [KoMa94], 6.9.2.

\medskip

{\bf 4.3.2. Definition.} {\it The structure of $\Cal{M}$--algebra on $L$ is a morphism of functors
$\bold{OpEnd}\,L\, \to H^*\Cal{M}$ compatible with glueing (see [KoMa94], 6.10.})

\medskip

The glueing itself and the compatibility  were introduced in [KoMa94]
in a somewhat  {\it ad hoc} manner. In fact glueing is the image of   {\it  ``grafting''} which must be included in
the list of basic morphisms of graphs, but then the minimal set of objects
upon which it is defined in order to deal with our (co)operad must include disjoint 
unions of stable graphs. For details, see [BorMa07], Sec.~1, and [KaWa17], Appendix A.

\medskip

{\bf 4.3.3. Proposition.} {\it The ``quantised action'' of $P(n)$ upon a quadratic algebra $Q$
(quantum linear space in the sense of [Ma88]) is
represented by the family of Kelly enrichments $P(n)^! \circ \underline{Hom}\, (Q^{\otimes n}, Q^{\otimes n})$
endowed with a family of elements similar to (4.6).
}
\medskip

{\bf Proof.} This directly follows from the Proposition 3.2 above.
I will supply some details that might not be immediately evident.

\medskip

Unfortunately, in the vast supply of examples of $P$--algebras, furnished
by quantum cohomology, I was unable to find nontrivial actions of $P$ upon {\it quadratic algebras} $A$
rather than upon {\it graded spaces} obtained by forgetting multiplication in $A$.

\smallskip

\medskip
{\bf 4.4. Cyclic $hyperCom$--algebras.} Generally, an operad can be characterised
by the category of algebras that it classifies.

\smallskip

The operad $P$ produces algebras endowed with infinitely many multilinear
operations satisfying infinitely many ``multicommutativity'' properties which I will briefly
recall below.

\smallskip

Let $L$ be a linear (super)space with symmetric even
non--degenerate scalar
product  $h$. 
\smallskip

An action of $P$ upon it induces upon $L$ the structure that we will
call here, following E.~Getzler, hypercommutative (or {\it hyperCom}) algebra.
In [Ma99], they were called $Comm_{\infty}$--algebras, but since then
the subscript $\infty$ acquired a standard connotation with cofibrant resolutions
(for operads) or homotopy lifts (for algebras).

\medskip

{\bf 4.4.1. Definition.}  {\it A structure of cyclic $hyperCom$--algebra on $(L,g)$
is a sequence of polylinear multiplications
$$
\circ_n:\ L^{\otimes n}\to L,\ \circ_n (\gamma_1\otimes\dots \otimes\gamma_n)=:
(\gamma_1,\dots ,\gamma_n),\ n\ge 2
$$
satisfying three axioms:

\medskip

(i) Commutativity = $\bold{S}_n$--symmetry;

\medskip

(ii) Cyclicity:\ $h((\gamma_1,\dots ,\gamma_n),\gamma_{n+1})$ is 
$\bold{S}_{n+1}$--symmetric;

\medskip

(iii) Associativity:\ for any $m\ge 0$, $\alpha ,\beta ,\gamma ,\delta_1,
\dots ,\delta_m$ 

$$
\sum_{\{ 1,\dots ,m\}=S_1\amalg S_2}
\pm ((\alpha ,\beta ,\delta_i\,|\,i\in S_1),\gamma ,\delta_j\,|\,j\in S_2)=
$$
$$
\sum_{\{ 1,\dots ,m\}=S_1\amalg S_2}
\pm (\alpha ,\delta_i\,|\,i\in S_1),\beta, \gamma ,\delta_j\,|\,j\in S_2))
$$
with usual signs from superalgebra.
\medskip

(iv) (Optional) identity Data and Axiom: $e\in L_{even}$ satisfying

\smallskip
$$
(e,\gamma_1,\dots ,\gamma_n)=\gamma_1\ \text{for}\ n=1;\ 0\ \text{for}\ n\ge 2.
$$
}
\medskip

This direct description of
cyclic $hyperCom$--algebras produces the same family of algebras
that was described in sec.~4.3.1 as $\Cal{M}$--algebras.

\medskip

{\bf 4.4.2. Examples.} 1) If $\circ_n=0$ for $n\ge 3$, we get the structure of commutative
algebra with invariant scalar product: $g(\alpha\beta ,\gamma )=
g(\alpha , \beta\gamma ).$

\bigskip

2) Associativity identities for $m=1$: 
$$
((\alpha ,\beta ),\gamma ,\delta )+((\alpha ,\beta ,\delta ),\gamma)=
((\alpha ,(\beta ,\gamma ,\delta ))+(\alpha ,\delta ,(\beta ,\gamma ))
$$

\smallskip

3) One of the earliest results of mathematical theory of quantum cohomology
established that for any smooth projective manifold 
(or a compact symplectic manifold) $V$, the superspace 
$$
(L,h):= (H^*(V),\ \text{Poincar\'e\ pairing})
$$
admits a canonical structure of cyclic $hyperCom$--algebra.

\medskip

{\bf 4.5. On self--reflexivity of the tree level quantum cohomology.} As I have already mentioned,
the idea to introduce a higher level (``quantised'') operadic action of $P$  upon its own components $\{P(n)\}$ was motivated by the problem which seems as yet far away from its solution: cf. [MaSm14]
and [MaSm13]. In the language of classical algebraic geometry, this problem consists
in calculation of Gromov--Witten invariants of genus zero of $\overline{M}_{0,n},\ n\ge 6$,
corresponding to those effective curve classes $\beta$ which lie ``to the wrong side'' of the  
anticanonical hyperplane, cf. [FarGi03].

\bigskip

\centerline{\bf References}

\medskip

%[BaWe85] M.~Barr, C.~Wells. {\it Toposes, Triples and Theories.} Springer, 1985.
 
%[3]   C.~Berger,  M.~Dubois--Violette, M.~Wambst. {\it Homogeneous
% algebras.} J. Algebra 261 (2003), no.~1, 172--185.
 
[BorMa07] D.~Borisov, Yu.~Manin. {\it Generalized 
operads and their inner cohomomorphisms.} Birkh\"auser Verlag, Progress in Math., vol. 265 (2007), 
247--308.

\smallskip

[BoyDr13] M.~Boyarchenko, V.~Drinfeld. {\it A duality formalism in the spirit 
 of Grothendieck and Verdier.} Quantum Topology, 4 (2013), 447--489.

% \smallskip
 
%[DavHe14]   L.~David, C.~Hertling. {\it Regular $F$--manifolds: initial conditions and Frobenius metrics.}
%arXiv:1411.4553v3

%\smallskip
%[6]   L.~David, I.~A.~B.~Strachan. {\it Dubrovin's duality for $F$--manifolds
%with eventual identities.} Advances in Math. vol.~206 (2011), 4031--4060.
\smallskip

 [DaSt04]   B.~Day, R.~Street. {\it Quantum categories, star autonomy, and quantum groupoids.}
Galois theory, Hopf algebras, and semiabelian categories, 187--225, Fields Inst. 
Commun., 43, Amer. Math. Soc., Providence, RI (2004), 187--225. 

\smallskip

[DeMi82]  P.~Deligne, J.~S.~Milne. {\it Tannakian categories.} In: Hodge cycles, motives, and Shimura varieties,
Springer Lecture Notes in Math, 900 (1982), 101--228.

\smallskip

[DoShVa12]   V.~Dotsenko, S.~Shadrin, B.~Vallette. {\it De Rham cohomology and
homotopy Frobenius manifolds.} Journ. Eur. Math. Soc., vol.~17, no.~3 (2015). arXiv:1203.5077.

\smallskip

[DoShVa15]   V.~Dotsenko, S.~Shadrin, B.~Vallette. {\it Noncommutative 
$\overline{\Cal{M}}_{0,n+1}.$} arXiv:1510.03261.

\smallskip

[Dr86] V.~Drinfeld. {\it Quantum groups.} Proc.~Int.~Congr.~Math., Berkeley (1986), vol.~1, 798--820. 

\smallskip

%[9]   B.~Dubrovin. {\it Geometry of 2D topological field theory.}
%In: Springer Lecture Notes in Math. 1620 (1996), 120--348.

%\smallskip

%[10]   B.~Dubrovin. {\it On almost duality for Frobenius manifolds.}
%Geometry, topology and math phys., Amer. Math. Soc. Translations, Ser. 1, 212 (2004), 75--132.
%arXiv: DG/0307374

% \smallskip

 [FadReTa87]  L.~D.~Faddeev, N.~Reshetikhin, L.~Takhtajan. {\it Quantization of Lie groups
 and Lie algebras.} Preprint LOMI (1987).
 
 \smallskip
 [FarGi03] G.~Farkas, A.~Gibney. {\it The Mori cones  of moduli spaces of
 pointed curves of small genus.}  TrAMS, vol. 355, No. 3 (2003), 1183--1199.
% \smallskip
 
% [GeKa95] E.~Getzler, M.~Kapranov. {\it Cyclic operads and cyclic homology.} In: Geometry, Topology,
 %and Physics for Raoul, ed.~by B.~Mazur, Internat.~Press, Cambridge MA (1995), 167--201.
 \smallskip
 
 [Fu84] W.~Fulton. {\it Intersection theory.} Springer Verlag, Berlin Heidelberg (1984), xi + 470 pp.
 
 \smallskip
 
 [GiKa94]    V.~Ginzburg, M.~Kapranov. {\it Koszul duality for operads.}
 Duke Math. J. 76 (1994), no. 1, 203--272. {\it Erratum:}
 Duke Math. J. 80 (1995), no. 1, 293.
 
% \smallskip

%[He02]  C.~Hertling. {\it Frobenius manifolds and moduli spaces for singularities.}
%Cambridge University Press, 2002.

%\smallskip

%[13]   C.~Hertling, Yu.~Manin. {\it Weak Frobenius manifolds.}
%Int. Math. Res. Notices, 6 (1999), 277--286. 
%arXiv:math.QA/9810132

%\smallskip

%[HeMaTe]   C.~Hertling, Yu.~Manin, C.~Teleman. {\it An update on semisimple quantum cohomology
%and $F$--manifolds.} Proc.~Steklov  Math.~Inst., vol.~264 (2009), 62--69. arXiv:math.AG/0803.2769

\smallskip

[Ji85] M.~Jimbo. {\it A $q$--difference analogue of $U(g)$ and the Yang--Baxter
equation.} Lett.~Math.~Phys. 10(1985), 247--252.

\smallskip

[KaWa17] R.~Kaufmann, B.~Ward. {\it Feynman categories.} Ast\'erisque 387, 2017. arXiv:1312.1269,
131 pp.
\smallskip
[Ke82]  G.~M.~Kelly. {\it Basic concepts of the enriched category theory.} Cambridge UP (1982).
Revised online version 
http//www.tac.mta.ca/tac/reprints/articles/10/tr10.pdf
\smallskip

[KoMa94] M.~Kontsevich, Yu.~Manin. {\it Gromov--Witten classes, quantum cohomology, and enumerative 
geometry}.  Comm. Math. Phys.,
164:3 (1994), 525--562.
%\smallskip

%[LoMa00] A.~Losev, Yu.~Manin. {\it New moduli spaces of pointed curves and pencils of flat connections.} Fulton's %Festschrift,
%Michigan Journ. of Math., 48 (2000), 443--472. arXiv:math.AG/0001003

\smallskip
 
 [LoVa12]  J.-L.~Loday, B.~Vallette. {\it Algebraic operads.} Springer (2012), xxiv+634 pp.
 \smallskip
[Ma87] Yu.~Manin. {\it Some remarks on Koszul algebras and
quantum groups.} Ann. Inst. Fourier, Tome XXXVII, f.~4 (1987), 191--205.

\smallskip
[Ma88] Yu.~Manin. {\it Quantum groups and non--commutative geometry.}
Publ. de CRM, Universit\'e de Montr\'eal (1988), 91 pp.
\smallskip
[Ma91] Yu.~Manin. {\it Topics in noncommutative geometry.} 
Princeton University Press (1991), 163 pp.

\smallskip

[Ma99] Yu.~Manin. {\it Frobenius manifolds, quantum cohomology, and moduli
spaces.} AMS Colloquium Publications, Vol.~47 (1999), xiii+303 pp.

\smallskip

[Ma17]   Yu.~Manin.  {\it Grothendieck--Verdier duality patterns in quantum algebra.}
 Izvestiya: Mathematics, vol. 81, No.~4, 2017. DOI: 10.4213/im8620 arXiv: 1701.01261. 14 pp.

\smallskip
[MaSm13]  Yu.~Manin, M.~Smirnov.  {\it On the derived category of $\overline{M}_{0,n}$.}
Izvestiya of Russian Ac. Sci., vol.~77, No 3 (2013), 93--108.
Preprint arXiv:1201.0265

\smallskip

[MaSm14]  Yu.~Manin,  M.~Smirnov.  {\it  Towards motivic quantum cohomology of $\overline{M}_{0,S}$.} 
Proc. of the Edinburg Math. Soc., Vol. 57 (ser. II), no 1, 2014,
pp. 201--230. Preprint arXiv:1107.4915

\smallskip

[MarShSt02] M.~Markl, St.~Shnider, J.~Stasheff. {\it Operads in algebra, topology and physics.}
 Math.~Surveys and Monographs, Vol.~96 (2002), x+349 pp.

% \smallskip

%[18]  S.~Merkulov. {\it Operads, deformation theory and $F$--manifolds.} In: Frobenius Manifolds. Quantum
%Cohomology and Singularities. Eds. C.~Hertling, M.~Marcolli. Aspects of Math. Vol.~ E36 (2004), 213--251.
%arXiv:math/0210478

\smallskip

[Va08]  B.~Vallette. {\it Manin products, Koszul duality, Loday algebras
and Deligne conjecture.} J. Reine Angew. Math. 620 (2008), 105--164.
arXiv:math/0609002

\bigskip

{\bf Max--Planck--Institute for Mathematics, 

Vivatsgasse 7, Bonn 53111, Germany.}

\smallskip

manin\@mpim-bonn.mpg.de

\enddocument